\documentclass[12pt]{amsart}
\usepackage{epsfig,color}
\usepackage{amsfonts, amsthm, amsmath}

\theoremstyle{plain}

\begin{document}

\title[Seven-game series vs. five-game series]{Are seven-game baseball playoffs fairer than five-game series when home-field advantage is considered?}
\author{Brian Dean}
\address{Department of Mathematics and Computer Science\\
    Salisbury University\\
    Salisbury, MD  21801}
\email{bjdean@salisbury.edu}
\date{}

\begin{abstract}
Conventional wisdom in baseball circles holds that a seven-game playoff series is fairer than a five-game series.  In an earlier paper, E. Lee May, Jr. showed that, treating each game as an independent event, a seven-game series is not significantly fairer.  In this paper, we take a different approach, taking home-field advantage into account.  That is, we consider a given series to consist of two disjoint sets of independent events---the home games and the road games.  We will take the probability of winning a given road game to be different from the probability of winning a given home game.  Our analysis again shows that a seven-game series is not significantly fairer.
\end{abstract}

\maketitle

\section{Introduction}\label{intro}
It is often said in baseball that a seven-game playoff series is fairer than a five-game series.  The argument is that, in a five-game series, the team without home-field advantage need only win its two home games, and take just one out of three on the road, in order to win the series.  On the other hand, to win a seven-game series, the team would have to either win all three home games and one of four on the road, or win at least two out of four on the road.  

Analyzing this question is a useful exercise in mathematical modeling and probability.  In \cite{Ma}, E. Lee May, Jr. showed that a seven-game series is not significantly fairer.  (By \textit{significantly fairer}, we mean that there is at least a four percent greater probability of winning the seven-game series than winning the five-game series.)  May approached the problem as follows: he let $p$ be the probability that the better team would win a given game in the series, and treated each game equally as an independent event without regard to where the game was being played.

In this paper, we will examine the same problem while attempting to account for home-field advantage.  From now on, $p$ will represent the probability that the team with home-field advantage in the series will win a given home game.  The probability that that team will win a given road game will be $rp$, where $r$, the \textit{road multiplier}, will be discussed in Section~\ref{roadmultiplier}.  Each home game will be treated as an independent event, and each road game will be treated as an independent event.

Since May approached the problem from the point of view of the better team, he necessarily had $p\in [0.5,1]$.  In this paper, where we approach the problem from the point of view of the team with home-field advantage, that will still be the case most of the time---in the Division Series and League Championship Series, home-field advantage goes to the better team.  However, in the World Series, this is not always the case.  Home-field advantage in the World Series alternated between the American and National Leagues through 2002; since 2003, it has been given to the champion of the league which had won that year's All-Star Game.  Still, in most cases, if a team is good enough to reach the World Series, then the probability that it will win a given home game is still likely to be at least 0.5, regardless of the opposition.  Nevertheless, it is possible that $p$ could be below 0.5, so we will only require $p\in [0,1]$.  Practically speaking, it seems unlikely that $p$ would ever be below, say, 0.4, but we will not require that to be the case.

\section{The Road Multiplier}\label{roadmultiplier}
As discussed in the Introduction, we will take the probability that the team with home-field advantage will win a given road game to be $rp$, where $r$ is a fixed number which we will call the \textit{road multiplier}.  For an individual team, the road multiplier is obtained by dividing the team's road winning percentage by its home winning percentage, i.e.,
$$\mbox{road multiplier}\,\, =\frac{\frac{RW}{RW+RL}}{\frac{HW}{HW+HL}}$$
where $RW$, $RL$, $HW$, and $HL$, are the number of the team's road wins, road losses, home wins, and home losses, respectively, in that season.  

Our value $r$ will be the average of the road multipliers of the 96 teams which have made the playoffs in the wildcard era (1995-2006).  This ends up giving us (to 9 decimal places)
$$r=0.894762228,$$
that is, we will consider the team with home-field advantage to be about 89.5 percent as likely to win a given road game as they are to win a given home game.

We will not list the results for all 96 teams here.  However, we will make a few comments.

The five highest and five lowest road multipliers of the 96 are as follows:

\bigskip

\begin{tabular}{lccc}
Team & Home Record & Road Record & Road Multiplier \\
\hline
2001 Braves & 40-41 & 48-33 & 1.2 \\
1997 Orioles & 46-35 & 52-29 & 1.130434783 \\
2001 Astros & 44-37 & 49-32 & 1.113636364 \\
2005 White Sox & 47-34 & 52-29 & 1.106382979 \\
2006 Tigers & 46-35 & 49-32 & 1.065217391 \\
(tie) 2000 White Sox & 46-35 & 49-32 & 1.065217391 \\
\hline
2000 Mets & 55-26 & 39-42 & 0.709090909 \\
2005 Braves & 53-28 & 37-44 & 0.698113208 \\
2006 Cardinals & 49-31 & 34-47 & 0.685311162 \\
2003 Athletics & 57-24 & 39-42 & 0.684210526 \\
2005 Astros & 53-28 & 36-45 & 0.679245283
\end{tabular}

\bigskip

Of the 96 teams, 23 of them had road multipliers of 1 or higher (meaning that about a quarter of the teams did at least as well on the road as they did at home), while 12 of the teams had road multipliers of 0.75 or below.  12 of the 16 highest road multipliers belong to American League teams, while 11 of the 16 lowest road multipliers belong to National League teams.  The road multipliers for the 12 World Series champions of the wildcard era, from highest to lowest, are as follows:

\bigskip

\begin{tabular}{lccc}
Team & Home Record & Road Record & Road Multiplier \\
\hline
2005 White Sox & 47-34 & 52-29 & 1.106382979 \\
1995 Braves & 44-28 & 46-26 & 1.045454545 \\
1999 Yankees & 48-33 & 50-31 & 1.041666667 \\
2000 Yankees & 44-35 & 43-39 & 0.941518847 \\
2001 Diamondbacks & 48-33 & 44-37 & 0.916666667 \\
1996 Yankees & 49-31 & 43-39 & 0.856147337 \\
1998 Yankees & 62-19 & 52-29 & 0.838709677 \\
2002 Angels & 54-27 & 45-36 & 0.833333333 \\
2004 Red Sox & 55-26 & 43-38 & 0.781818182 \\
1997 Marlins & 52-29 & 40-41 & 0.769230769 \\
2003 Marlins & 53-28 & 38-43 & 0.716981131 \\
2006 Cardinals & 49-31 & 34-47 & 0.685311162
\end{tabular}

\bigskip

\section{Comparing Three-Game Series and Five-Game Series}\label{threeversusfive}

Before comparing seven-game series and five-game series, we will first look at five-game series versus three-game series, as that case is a bit easier to dive right into.  Throughout the next two sections, we will use the following notation: we will use capital letters (W and L) to denote games in which the team with home-field advantage wins and loses at home, and lowercase letters (w and l) to denote games in which that team wins and loses on the road.  Thus, each instance of W will have probability $p$, each L will have probability $1-p$, each w will have probability $rp$, and each l will have probability $1-rp$.

\subsection{Three-game series}

There have never been three-game playoff series in baseball, except to break ties (most notably the playoff between the New York Giants and Brooklyn Dodgers following the 1951 season).  However, if there were, they would likely be in one of two formats---either a 1-1-1 format (in which the team with home-field advantage plays games one and three at home, and game two on the road), or a 1-2 format (in which they play game one on the road and games two and three at home).  

The scenarios for that team to win the series, in a 1-1-1 format, are as follows.

\bigskip

\begin{tabular}{lc}
Scenario & Probability \\
\hline
Ww & $p(rp)$ \\
WlW & $p^2(1-rp)$ \\
LwW & $p(rp)(1-p)$
\end{tabular}

\bigskip
Adding these probabilities, we see that the total probability that the team with home-field advantage will win the series, in a 1-1-1 format, is $(2r+1)p^2-2rp^3$.

The following are the corresponding scenarios if the series were played in a 1-2 format.

\bigskip

\begin{tabular}{lc}
Scenario & Probability \\
\hline
wW & $p(rp)$ \\
wLW & $p(rp)(1-p)$ \\
lWW & $p^2(1-rp)$
\end{tabular}

\bigskip
Again, the total probability of victory in this format is $(2r+1)p^2-2rp^3$.  So, the probability that the team with home-field advantage will win a three-game series is the same in either format.

\subsection{Five-game series}
Major League Baseball employed five-game playoff series for the League Championship Series from 1969-1984.  (Prior to 1969, the playoffs consisted solely of the teams with the best records in each league meeting in the World Series.)  Since 1985, the League Championship Series have been in a best-of-seven format.  However, five-game series returned with the advent of the wildcard system; since 1995, each league has had two five-game Division Series, with the winners advancing to the seven-game League Championship Series.

Two formats for best-of-five series have been used over the years: a 2-3 format (in which the team with home-field advantage plays the first two games on the road and the final three games at home), and a 2-2-1 format (in which that team plays games one, two, and five at home, and games three and four on the road).  We will examine each format separately; as with the two formats for three-game series, we will see that the probability that the team with home-field advantage will win the series is independent of the format.

First, we examine the scenarios in which the team with home-field advantage will win the series, if the series is in a 2-3 format.

\bigskip
\begin{tabular}{lc}
Scenario & Probability \\
\hline
wwW & $p(rp)^2$ \\
lwWW & $p^2(rp)(1-rp)$ \\
wlWW & $p^2(rp)(1-rp)$ \\
wwLW & $p(rp)^2(1-p)$ \\
llWWW & $p^3(1-rp)^2$ \\
lwLWW & $p^2(rp)(1-p)(1-rp)$ \\
lwWLW & $p^2(rp)(1-p)(1-rp)$ \\
wlLWW & $p^2(rp)(1-p)(1-rp)$ \\
wlWLW & $p^2(rp)(1-p)(1-rp)$ \\
wwLLW & $p(rp)^2(1-p)^2$
\end{tabular}

\bigskip

Summing these, we see that the total probability that the team with home-field advantage will win the series, in a 2-3 format, is
$$(3r^2+6r+1)p^3-(9r^2+6r)p^4+6r^2p^5$$

Next, we look at the corresponding scenarios for a 2-2-1 format.

\bigskip
\begin{tabular}{lc}
Scenario & Probability \\
\hline 
WWw & $p^2(rp)$ \\
LWww & $p(rp)^2(1-p)$ \\
WLww & $p(rp)^2(1-p)$ \\
WWlw & $p^2(rp)(1-rp)$ \\
LLwwW & $p(rp)^2(1-p)^2$ \\
LWlwW & $p^2(rp)(1-p)(1-rp)$ \\
LWwlW & $p^2(rp)(1-p)(1-rp)$ \\
WLlwW & $p^2(rp)(1-p)(1-rp)$ \\
WLwlW & $p^2(rp)(1-p)(1-rp)$ \\
WWllW & $p^3(1-rp)^2$
\end{tabular}

\bigskip
Again, if we add these, we see that the total probability of victory is
$$(3r^2+6r+1)p^3-(9r^2+6r)p^4+6r^2p^5$$
and so the probability that the team with home-field advantage will win a five-game series is the same in either format.

\subsection{Comparing the two}
To find the difference in probabilities in winning a five-game series and a three-game series, we just subtract the two: the probability of winning a five-game series, minus the probability of winning a three-game series, is the function 
$$f(p)=6r^2p^5-(9r^2+6r)p^4+(3r^2+8r+1)p^3-(2r+1)p^2,\;\;p\in [0,1]$$
We will find the extreme values of $f$ using the Extreme Value Theorem.

The derivative of $f$ is
$$f'(p)=30r^2p^4-(36r^2+24r)p^3+(9r^2+24r+3)p^2-(4r+2)p;$$
keeping in mind that $r=0.894762228$, the derivative is 0 for
$$p=0$$
$$p\approx 0.294269665$$
$$p\approx 0.756820873$$
and for a value of $p$ between 1 and 2 (as can be verified using the Intermediate Value Theorem).

Checking the values of $f$ at the critical points and the endpoints, we get
$$f(0)=0$$
$$f(0.294269665)\approx -0.056156576$$
$$f(0.756820873)\approx 0.047338476$$
$$f(1)=0$$
So, a five-game series is at most about 4.73\% fairer than a three-game series, and at worst about 5.62\% less fair.  However, as mentioned in the Introduction, it is extremely unlikely that $p$ would ever be as low as $0.294$.  If we look at the value of $f$ at a more realistic lower bound for $p$, we get
$$f(0.4)=-0.0431953192$$
and so a five-game series is about 4.32\% less fair than a three-game series for that value of $p$.  In summary, there does appear to be a significant difference between three-game and five-game series for certain values of $p$.

The value of $p$ in $[0,1]$ for which $f(p)=0$ is approximately 0.537783; for $p$ less than that, three-game series are fairer (or, put another way, five-game series are fairer for the team without home-field advantage), while for $p$ greater than that, five-game series are fairer.

\section{Comparing Five-Game Series and Seven-Game Series}\label{fiveversusseven}
We are now ready to examine the question of interest to us, comparing a five-game series and a seven-game series.  We have already shown that the probability that the team with home-field advantage will win a five-game series, regardless of format, is
$$(3r^2+6r+1)p^3-(9r^2+6r)p^4+6r^2p^5$$

\subsection{Seven-game series}
A seven-game series in baseball is played under a 2-3-2 format---the team with home-field advantage plays games one, two, six, and seven at home, and the middle three games on the road.  There are a total of 35 possible scenarios for victory, so we will not list each separately.  However, we will list each scenario lasting four, five, or six games.

\bigskip
\begin{tabular}{lc}
Scenario & Probability \\
\hline
WWww & $p^2(rp)^2$ \\
LWwww & $p(rp)^3(1-p)$ \\
WLwww & $p(rp)^3(1-p)$ \\
WWlww & $p^2(rp)^2(1-rp)$ \\
WWwlw & $p^2(rp)^2(1-rp)$ \\
LLwwwW & $p(rp)^3(1-p)^2$ \\
LWlwwW & $p^2(rp)^2(1-p)(1-rp)$ \\
LWwlwW & $p^2(rp)^2(1-p)(1-rp)$ \\
LWwwlW & $p^2(rp)^2(1-p)(1-rp)$ \\
WLlwwW & $p^2(rp)^2(1-p)(1-rp)$ \\
WLwlwW & $p^2(rp)^2(1-p)(1-rp)$ \\
WLwwlW & $p^2(rp)^2(1-p)(1-rp)$ \\
WWllwW & $p^3(rp)(1-rp)^2$ \\
WWlwlW & $p^3(rp)(1-rp)^2$ \\
WWwllW & $p^3(rp)(1-rp)^2$
\end{tabular}

\bigskip
There are a total of 20 scenarios for victory which last the full seven games.  Rather than list each one separately, we will just list the various combinations of W, L, w, and l, give the probability of each occurrence, and give the number of ways each scenario occurs.  For example, occurrences of the first type include LLlwwWW, LWwlwLW, and WLwwlLW.

\bigskip
\begin{tabular}{lcc}
Scenario & Probability & Occurrences \\
\hline
2 W, 2 w, 2 L, 1 l & $p^2(rp)^2(1-p)^2(1-rp)$ & 9 \\
3 W, 1 w, 1 L, 2 l & $p^3(rp)(1-p)(1-rp)^2$ & 9 \\
1 W, 3 w, 3 L, 0 l & $p(rp)^3(1-p)^3$ & 1 \\
4 W, 0 w, 0 L, 3 l & $p^4(1-rp)^3$ & 1
\end{tabular}

\bigskip
Adding together all of the probabilities for the 35 victory scenarios, we see that the total probability that the team with home-field advantage will win a seven-game series is
$$(4r^3+18r^2+12r+1)p^4-(24r^3+48r^2+12r)p^5+(40r^3+30r^2)p^6-20r^3p^7$$

\subsection{Comparing the two}
If we take the probability of winning a seven-game series, and subtract the probability of winning a five-game series, we get the function
\begin{eqnarray*}
s(p) &=& -20r^3p^7+(40r^3+30r^2)p^6-(24r^3+54r^2+12r)p^5 \\
  & & \,\, +(4r^3+27r^2+18r+1)p^4-(3r^2+6r+1)p^3
\end{eqnarray*}
where $p\in [0,1]$.

The derivative of this function is
\begin{eqnarray*}
s'(p) &=& -140r^3p^6 + (240r^3+180r^2)p^5-(120r^3+270r^2+60r)p^4 \\
  & & \,\, +(16r^3+108r^2+72r+4)p^3-(9r^2+18r+3)p^2
\end{eqnarray*}
Again using the fact that we are taking $r=0.894762228$, the derivative $s'$ is 0 for
$$p=0$$
$$p\approx 0.329786090$$
$$p\approx 0.723663130$$
and for a value of $p$ between 1 and 1.05, and a value of $p$ between 1.05 and 1.1.  (These last two can be verified using the Intermediate Value Theorem.)

Checking the values of $s$ at the critical points and the endpoints, we get
$$s(0)=0$$
$$s(0.329786090)\approx -0.038565024$$
$$s(0.723663130)\approx 0.034221072$$
$$s(1)=0$$
So, a seven-game series is at most about 3.42\% fairer than a five-game series, and at worst about 3.86\% less fair (and that occurs for a value of $p$ which is likely too small to occur in practice).  Therefore, there is no significant difference between a five-game series and a seven-game series.

The value of $p$ in $[0,1]$ for which $s(p)=0$ is approximately 0.533711; for $p$ less than that, five-game series are fairer (i.e., seven-game series are fairer for the team without home-field advantage), while for $p$ greater than that, seven-game series are fairer.

\section{Further Questions}\label{furtherquestions}

There are a few ways in which this model could be amended.  First, instead of finding a fixed value of $r$ for the road multiplier, we could keep $r$ as a variable (with appropriate upper and lower bounds for $r$), and then treat the functions $f$ and $s$ as functions of two variables.

Another approach would be to account for morale.  In \cite{Re}, S. Reske approaches the problem as May did in \cite{Ma}---that is, with $p$ representing the probability that the better team would win a given game, without regard to home-field advantage.  However, if the better team has a lead in the series, then its probability of winning the next game would be $p+a$, while if it trails in the series, then its probability of winning the next game would be $p-a$, where $a$ may be either positive or negative.  The idea is that, if the team leads the series, its increase in morale (and subsequent decrease in the other team's morale) could actually make it more likely to win the next game, and vice versa if it trails the series.  In that case, $a>0$.  The case $a<0$ would correspond to what happens if the team leads the series but then gets overconfident, making it less likely to win the next game.  With this approach, Reske again shows that there is no significant difference between a five-game series and a seven-game series.  This could be easily adapted to account for home-field advantage, with the fixed value of $r$ we used in this paper: if the team with home-field advantage leads the series, and the next game is at home, its probability of winning would be $p+a$, while if the next game were on the road, it would be $r(p+a)$; similarly if the team with home-field advantage trails the series, its probability of winning the next game would be $p-a$ if at home, and $r(p-a)$ if on the road.  This would again be a two-variable problem, with variables $p$ and $a$.  If we do not require $r$ to be fixed, then it would become a three-variable problem.

A final approach could be one of cumulative morale.  That is, if the team with home-field advantage leads the series by one game, then its probability of winning the next game would be $p+a$ or $r(p+a)$, if it leads the series by two games, its probability of winning the next game would be $p+2a$ or $r(p+2a)$, and so forth.  The idea here would be that, the further ahead the team is, the greater its morale would get (if $a>0$), or the more overconfident it would get (if $a<0$).


\begin{thebibliography}{99}

\bibitem{Ma} E.L. May, Jr., {\em Are Seven-Game Baseball Playoffs Fairer?}, Mathematics Teacher, {\bf\underline{85}}, (Oct. 1992), 528--531. 

\bibitem{Re} S. Reske, {\em Are Seven Game Baseball Series Still Fairer when Morale is Considered?}, preprint, http://www4.carthage.edu/faculty/ msnavely//sdlnew/vol1/baseball.pdf.

\end{thebibliography}
\end{document}